\documentclass[a4paper,12pt]{article}

\usepackage{amsmath}
\usepackage{amssymb}
\usepackage{amsthm}
\usepackage{mathtools}
\usepackage{mathrsfs}
\usepackage{graphicx}
\usepackage{subcaption}
\usepackage{color}
\usepackage{multirow}
\usepackage[hmargin=2.5cm,vmargin=2.5cm]{geometry}

\newtheorem{proposition}{Proposition}
\newtheorem{remark}{Remark}

\author{H. Ghane\footnote{Bandar Anzali Branch, Islamic Azad University, Bandar Anzali, Iran (\texttt{h.ghane.sasansaraei@rug.nl})} \footnotemark[2] \and A.E. Sterk\footnote{Bernoulli Institute, University of Groningen, Groningen, The Netherlands (\texttt{a.e.sterk@rug.nl}, \texttt{h.waalkens@rug.nl})} \and H. Waalkens\footnotemark[2]}
\title{Dynamical analysis of a chaos generator}
\date{\today}

\begin{document}

\maketitle

\begin{abstract}
{Investigating the possibility of applying techniques from linear systems theory to the setting of nonlinear systems has been the focus of many papers. The pseudo linear form representation of nonlinear dynamical systems has led to the concept of nonlinear eigenvalues and nonlinear eigenvectors. When the nonlinear eigenvectors do not depend on the state vector of the system, then the nonlinear eigenvalues determine the global qualitative behaviour of a nonlinear system throughout the state space. The aim of this paper is to use this fact to construct a nonlinear dynamical system of which the trajectories of the system show continual stretching  and folding. We first prove that the system is globally bounded. Next, we analyse the system numerically by studying bifurcations of equilibria and periodic orbits. Chaos arises due to a period doubling cascade of periodic attractors. Chaotic attractors are presumably of H\'enon-like type, which means that they are the closure of the unstable manifold of a saddle periodic orbit. We also show how pseudo linear forms can be used to control the chaotic system and to synchronize two identical chaotic systems.}
\end{abstract}

\newpage


\section{Introduction}

The analysis of nonlinear systems is a wide field of research with many applications and techniques. 
One  main approach to the analysis and control of nonlinear systems
 consists  of transferring results from linear systems theory.  
The best known example is the Poincar{\'e} linearization near an equilibrium point where for a hyperbolic equilibrium point, the linear dynamics associated with the Jacobian matrix of the vector field is by the Hartman-Grobman Theorem conjugate to the nonlinear dynamics near the equilibrium point, see, e.g.,  Cheng {\it et al.}\ (2010). 
As another linearization scheme we mention  feedback linearization which amounts to designing a feedback control along with some change of coordinate which transforms the closed loop nonlinear system into a linear system, see Baillieul \& Willems (1999).

However, one of the most effective applications of linear systems theory in nonlinear systems is the State-Dependent Riccati Equation (SDRE) strategy in nonlinear optimal control theory, see \c{C}imen (2008). This approach requires a representation of the nonlinear dynamics into a linear form with a state dependent system matrix. In doing so, this matrix valued function fully captures the nonlinearities of the system, which provides the designer a very effective method of making a good and yet systematic trade-off between state error and input effort via a state dependent linear quadratic formulation. A SDRE, of which the coefficients vary across state space, is then solved to give a suboptimal control law. The SDRE approach in nonlinear optimal control design relies on the pseudo linear (PL) representation of a nonlinear dynamical system. Indeed, the closed loop system obtained by this optimal controller is still in a PL form. However, the stability analysis of the closed loop system by exerting the resulting optimal control is still a problematic challenge, and has attracted several studies during the last years. As noted by Cloutier (1997), the number of successful applications of the SDRE approach in the design of nonlinear optimal controllers outpaced the available theoretical results.

Investigating the possibility of using the PL form representation in the stability analysis of nonlinear systems has been the effort of several works, e.g.\ Banks \& Mhana (1992), Tsiotras {\it et al.}\ (1996), Banks \& Mhana (1996), Langson \& Alleyne (2002) and Muhammad \& Van Der Woude (2009). The key focus of these works was the stabilizability of a PL form by exerting the state dependent control obtained via the SDRE approach. Recently, Ghane \& Menhaj (2015) introduced a theorem providing a sufficient condition of a PL system for correct stability analysis based on its state dependent eigenvalues and eigenvectors. Although the PL representation was originally introduced for the systematic design of a nonlinear optimal controller through the SDRE approach, Ghane \& Menhaj (2015) have shown that besides stability analysis, the PL form can also provide a useful tool for the global qualitative analysis of nonlinear dynamical systems when the nonlinear eigenvectors obtained from this PL representation are state independent.

The ability of determining the qualitative behaviour of a nonlinear dynamical system by means of eigenvalues and -vectors obtained from a PL form is also attractive for fields beyond control engineering applications, such as dynamical systems (see Ghane \& Menahj (2014)). The wide spread applications of chaotic systems in practical applications like image watermarking (see Wang {\it et al.}\ (2015)), chaotic communication (see \c{C}i\c{c}ek {\it et al.}\ (2016) and Zhou {\it et al.}\ (2014)), robotics (see Zang {\it et al.}\ (2016)),  have motivated us to apply this qualitative approach to generate a class of chaotic systems. In this paper, we apply the PL representation in the interesting field of chaos generation which may be potentially useful for the engineering applications mentioned above.

The aim of this paper is to synthesize the basic qualitative characteristics of a chaotic behaviour. With the help of nonlinear eigenvalues as the qualitative indicator of the behaviour of a system in PL form, we introduce a system with a specific type of locally unstable and globally bounded trajectories. The system that we construct in this way has equilibria that become unstable through a Hopf bifurcation. The resulting periodic orbits bifurcate further through a period doubling cascade which leads to chaotic attractors. The latter are presumably of H\'enon-like type which means that they are the closure of the unstable manifold of a saddle periodic orbit. In addition, we show the application of nonlinear eigenvalues in nonlinear control and synchronization of the chaotic system.

The rest of the paper is organized as follows. In section 2, the PL representation of nonlinear systems is briefly introduced and its ability for qualitative analysis of nonlinear dynamical systems is discussed. Section 3 is devoted to use the PL form to generate a chaotic system and to prove the global boundedness of the trajectories. The dynamical analysis of the obtained chaotic system is presented in section 4. In section 5, an eigen-structure based analysis is used to design a control law for the chaotic system and to perform an identical synchronization of two chaotic systems. Finally, some concluding remarks are presented in section 6.


\section{Pseudo linear systems: a brief review}

An autonomous nonlinear system is described by a system of nonlinear ordinary differential equations which do not explicitly depend on the independent variable. The general form of such a system is given by
\begin{equation}
	\label{eq:general}
	\bf{\dot x}(t) = \bf{f}(\bf{x}(t)),
\end{equation}
where $\bf{x}$ takes values in $n$-dimensional Euclidean space and the independent variable $t$ is usually time. Now assume that $\bf{f} : \mathbb{R}^n \to \mathbb{R}^n$ is sufficiently smooth and that $\bf{f}(\bf{0}) = \bf{0}$. Inspired by the theory of linear systems we can then transform an autonomous system \eqref{eq:general} to the form
\begin{equation}
	\label{eq:generalPL}
	{\bf{\dot x}}(t) = A({{\bf{x}}(t)}){\bf{x}}(t)
\end{equation}
where $A : \mathbb{R}^n \to \mathbb{R}^{n \times n}$. This form is called pseudo linear (PL) and it was originally introduced by Banks \& Mhana (1992) to cope with the difficulty of designing optimal control laws for nonlinear systems.

Using the PL form \eqref{eq:generalPL}, it is possible to extend the concept of eigenvalues and eigenvectors to the setting of nonlinear systems. The nonlinear eigenvalue (NEValue) and its corresponding nonlinear eigenvector (NEVector) are defined as the functions  $\lambda : \mathbb{R}^n \to \mathbb{C}$ and ${\bf v} : \mathbb{R}^n \to \mathbb{C}^n$, respectively, that satisfy the equation
\begin{equation}
	A(\bf{x})\bf{v}(\bf{x}) = \lambda (\bf{x})\bf{v}(\bf{x}).
\end{equation}
Equivalently, the nonlinear eigenvalues can also be obtained as the solution of the characteristic equation
\begin{equation}
	\det ({A({\bf{x}}) - \lambda ({\bf{x}}){I_n}}) = 0.
\end{equation}

Based on these generalized concepts, the following remarks and proposition are presented. Proofs and more explanations can be found in Ghane \& Menhaj (2014) and (2015). By means of these results we can study the qualitative behaviour of nonlinear dynamical systems. The qualitative analysis of nonlinear systems based on PL forms mainly uses the following observation from linear systems theory.

\begin{remark}
The qualitative behaviour of nonlinear systems is determined by:
\begin{enumerate}
    \item The sign of the real part of the NEValues;
   	\item The realness or complexness of the NEValues.
\end{enumerate}
The first condition determines the stability properties of the dynamics and the second condition determines the spiraling or exponential nature.
\end{remark}

\begin{remark}
Consider a nonlinear dynamical system $\bf{\dot x} = \bf{f}\left( \bf{x} \right)$, where $\bf{f} : \mathbb{R}^n \to \mathbb{R}^n$ satisfies $\bf{f}(\bf{0}) = \bf{0}$. Among the infinite distinct possible PL forms of this system, only the unique PL form which has state independent (SI) NEVectors must be used in eigen-structure based analysis, because only a PL form with SI NEVectors is guaranteed to yield correct qualitative results through its NEValues analysis.
\label{Pro1}	
\end{remark}

\begin{proposition}
For a nonlinear system $\bf{\dot x} = \bf{f}(\bf{x})$ where $\bf{f} : \mathbb{R}^n \to \mathbb{R}^n$ satisfies $\bf{f}(\bf{0}) = \bf{0}$, a sufficient condition for global asymptotic stability of the origin is that the system has a PL form representation that satisfies the following conditions:
\begin{enumerate}
	\item $\mathop{\rm Re}\nolimits \left\{ \lambda_i(\bf{x}) \right\} < 0$ for all $\bf{x} \in \mathbb{R}^n$ and $i=1,\dots,n$;
	\item For every NEValue the algebraic and geometric multiplicities are equal;
	\item All NEVectors of the matrix $A\left( {\bf{x}} \right)$ are state independent.
\end{enumerate}
\end{proposition}

\begin{remark}
All of the aforementioned results are applicable to the special class of nonlinear systems of the form \eqref{eq:generalPL} in which $A({\bf{x}}) = \mathop{diag}(D_1({\bf{x}}),\dots,D_p({\bf{x}}))$ is a block diagonal matrix and where each block is of the form
\[
D_i({\bf{x}}) = \begin{bmatrix} g_i({\bf{x}}) \end{bmatrix}
\quad\text{or}\quad
D_i({\bf{x}}) = \begin{bmatrix} g_i({\bf{x}}) & -\omega_i \\ \omega_i & g_i({\bf{x}}) \end{bmatrix},
\]
where $g_i : \mathbb{R}^n \to \mathbb{R}$ and $\omega_i > 0$ for $i=1,\dots,p$. Each $1\times 1$ block gives rise to a real NEValue and each $2 \times 2$ block gives rise to a complex NEValue. It is straightforward to verify that the corresponding NEVectors are state independent.
\end{remark}


\section{Chaos generation}

Qualitatively speaking, 
the occurrence of chaotic behaviour is usually related to the interplay between local instability and global boundedness of trajectories (see Sch\"oll \& Schuster (2008)). The local instability is responsible for the exponential divergence of nearby trajectories, whilst the global boundedness folds trajectories within the finite volume of the system's phase space. The combination of these two mechanisms can result in high sensitivity of the system trajectories to the initial conditions. In this paper we use NEValues of a PL form as indicators of a system's qualitative behaviour and to construct a dynamical system that shows chaotic behaviour without the need for exhaustive tuning of parameter values.


\subsection{Constructing a candidate chaotic system}

The Poincar\'e-Bendixson Theorem implies that the dynamical behaviour a system of the form \eqref{eq:general} with $n=2$ cannot be chaotic, see Guckenheimer \& Holmes (1983). Hence, the minimum dimension of a chaotic system is $n=3$. We first concentrate on finding nonlinear functions ${\lambda _i}({\bf{x}})$, where $i = 1,2,3$, to generate the continual stretching and folding property in the dynamics of the system. By the approach proposed in the previous section, it is possible to produce such a behaviour with a proper selection of NEValues. As a result, this approach may lead to different choices of NEValues which satisfy the desired qualitative behaviour; one of these choices is the following one:
\begin{equation}
	\begin{split}
	\lambda_{1,2}(\bf{x}) & = (x_3^2 - h^2) \pm j\omega,  \\
	\lambda_{3}(\bf{x}) & = r^2 - a x_1^2 - b x_2^2 - c x_3^2,
	\end{split}
\end{equation}
in which $a, b, c, \omega, h, r > 0$ are fixed parameters and $j^2=-1$.

Applying Remark 2.1, these NEValues give rise to the following nonlinear system:
\begin{equation}
\label{eq:ghane}
	\begin{split}
	\dot{x}_1 & = ({x_3}^2 - {h^2}){x_1} - \omega {x_2},  \\
	\dot{x}_2 & = \omega {x_1} + ({x_3}^2 - {h^2}){x_2},  \\
	\dot{x}_3 & = ({r^2} - a{x_1}^2 - b{x_2}^2- c{x_3}^2){x_3}.
	\end{split}
\end{equation}
The NEVectors of the system \eqref{eq:ghane}, which are simply given by
\[
	\bf{v}_1 = \begin{bmatrix} 1 & -j & 0 \end{bmatrix}^\top,	\quad
	\bf{v}_2 = \begin{bmatrix} j &  1 & 0 \end{bmatrix}^\top,	\quad
	\bf{v}_3 = \begin{bmatrix} 0 &  0 & 1 \end{bmatrix}^\top,
\]
satisfy the condition of Remark~\ref{Pro1}. Therefore, the chosen NEValues guarantee that the system will exhibit the continual stretching and folding that is characteristic of a chaotic system. 
For further analysis, it is convenient to set $x_1 = \rho\cos\theta$ and $x_2 = \rho\sin\theta$ by which the system \eqref{eq:ghane} can be rewritten in terms of cylindrical coordinates:
\begin{equation}
	\label{eq:ghane-cylinder}
	\begin{split}
	\dot{\rho}_{\phantom{3}}   & = (x_3^2 - h^2)\rho, \\
	\dot{\theta}_{\phantom{3}}  & = \omega, \\
	\dot{x}_3    & = (r^2 - a\rho^2\cos^2\theta - b\rho^2\sin^2\theta - cx_3^2)x_3.
	\end{split}
\end{equation}

\begin{figure}
\raisebox{4cm}{a)}\includegraphics[width=0.29\linewidth] {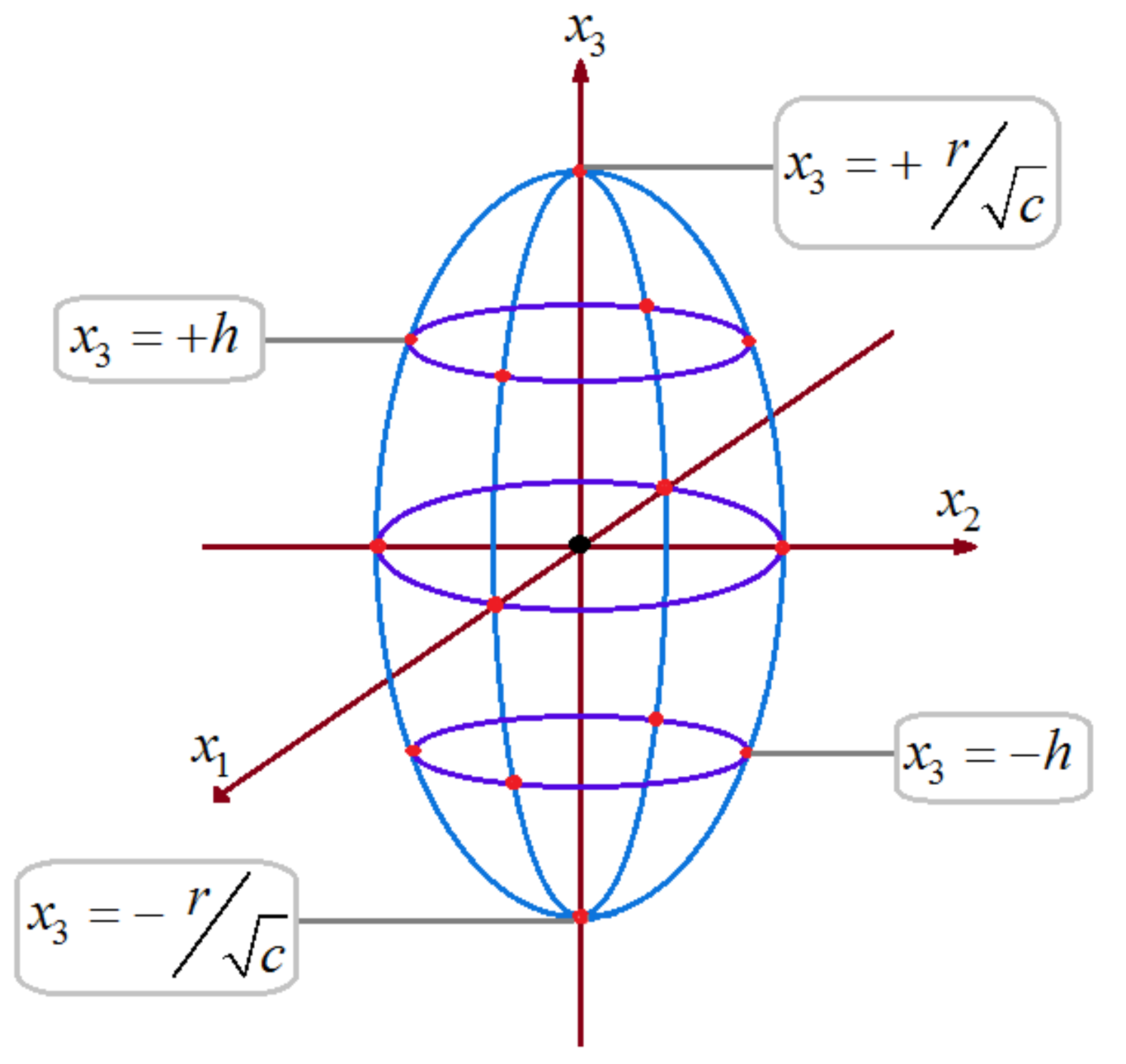}
\raisebox{4cm}{b)}\includegraphics[width=0.29\linewidth]{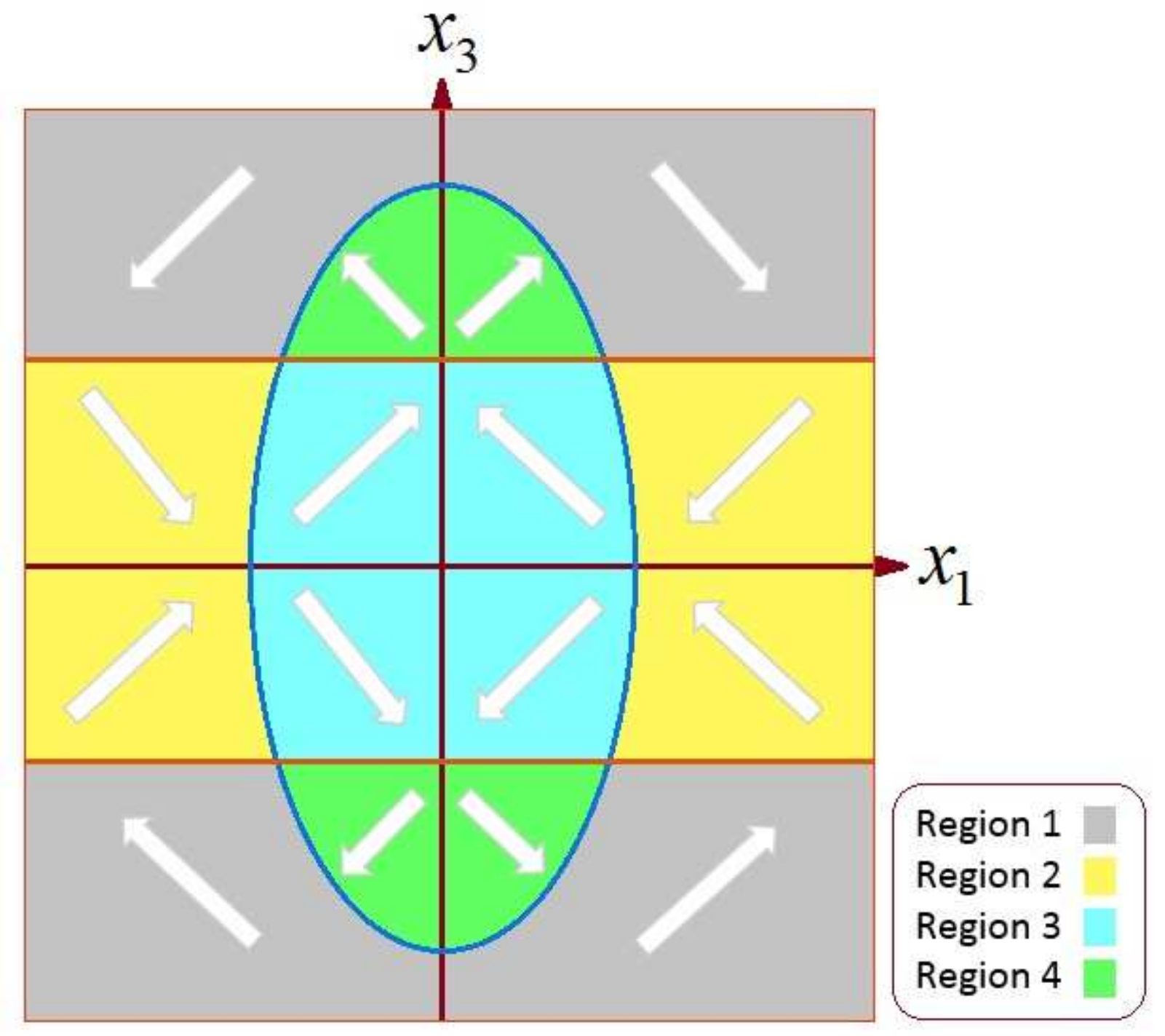}
\raisebox{4cm}{c)}\includegraphics[width=0.29\linewidth]{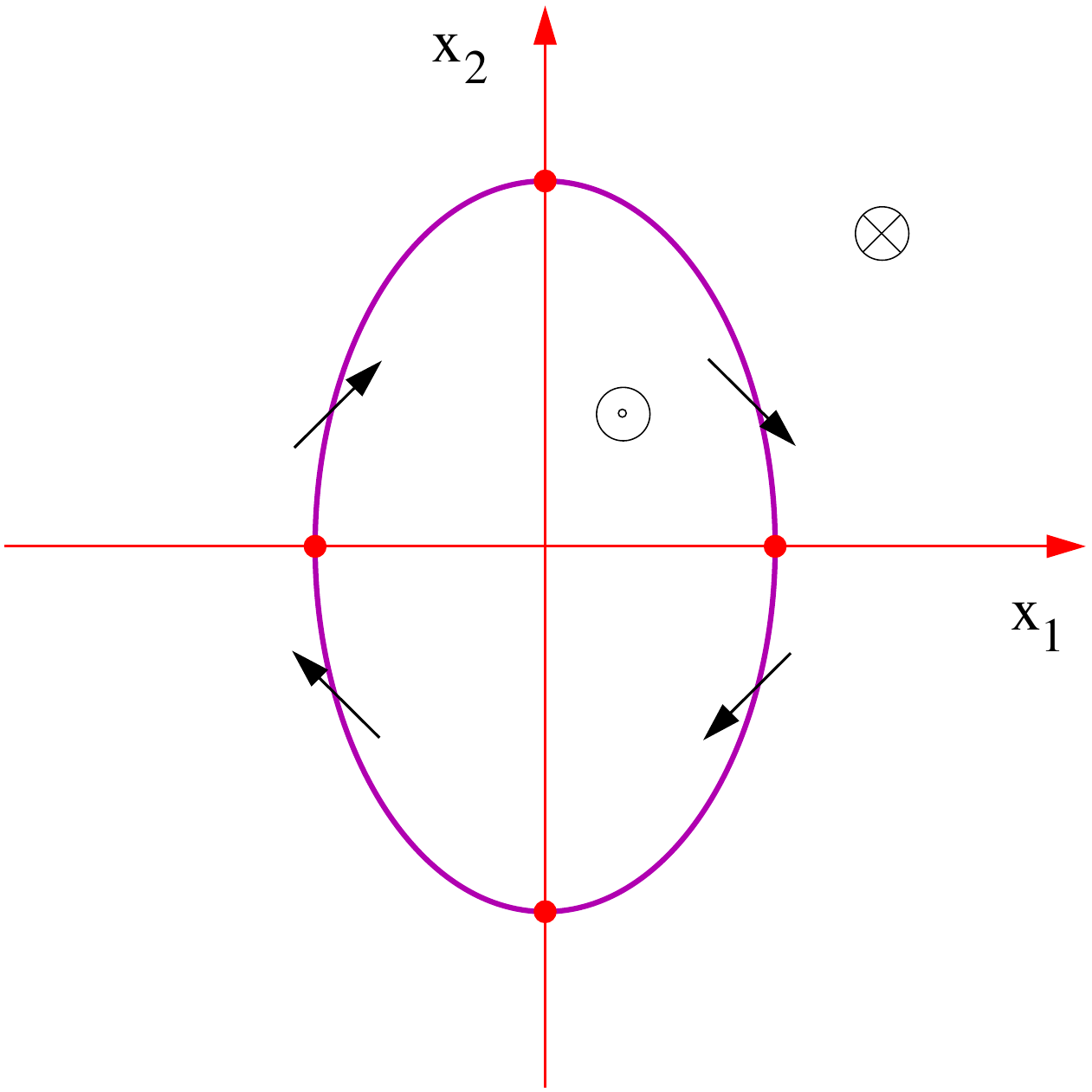}
\caption{Graphical illustration of the qualitative behaviour of system \eqref{eq:ghane}.
		On the ellipsoid $a{x_1}^2 + b{x_2}^2 + c{x_3}^2 = r^2$ and on the plane $x_3=0$ the $x_3$-component of the vector field is vanishing. On the planes $x_3=\pm h$, 
		$\dot{\rho}=0$. 
a) Nullclines $\dot{x}_3=0$ and $\dot{\rho}=0$ in the three-dimensional state space.
b) Section $x_2=0$ of the state space. The arrows indicate direction of the vector field in terms of the signs of $\dot{x}_3$ and $\dot{\rho}$ in the regions 1, 2, 3 and 4.
c) Section $x_3=0$ with the intersection of the nullclines along  the ellipse  $a{x_1}^2 + b{x_2}^2 + c{h}^2 = r^2$. In the region enclosed by  the ellipse the $x_3$-component of the vector field is positive. Outside it is negative. The arrows   indicate the direction of the vector field  on the ellipse.
}
\label{fig:Qualitative_analysis}
\end{figure}

As illustrated in Figure~\ref{fig:Qualitative_analysis}
the $\rho$ nullclines given by the horizontal planes $x_3=\pm h$ and the $x_3$ nullclines given by the ellipsoid $a{x_1}^2 + b{x_2}^2 + c{x_3}^2 = r^2$ and the horizontal plane 
$x_3=0$ divide the state space into four regions with different signs of the real parts of the NEValues 
that give rise to the different qualitative behaviours described in Table~\ref{tab:Qualitative}. Note that we consider the regions to be open, i.e. to not contain their boundaries. 
In region~2 and region~4, the real parts of all NEValues are positive and negative, respectively. Thus, in these regions the trajectories of the system are repelled from and attracted to the origin, respectively. Inline  with our approach for chaos generation, the existence of these two types of behaviours besides the regions 1 and 3 with saddle behaviour is necessary to ensure both the stretching and the folding of system trajectories. A proper arrangement of these regions then guarantees that the trajectories of the system remain bounded which is proved in detail on the next subsection. This arrangement is assured in system \eqref{eq:ghane} by the assistance of the regions 1 and 3.

Observe that replacing $x_3$ by $-x_3$ in \eqref{eq:ghane} yields the same equations. 
As a consequence the plane $x_3=0$ is invariant under the flow.
Therefore, it suffices to discuss the dynamics for $x_3 > 0$ and this is what we do for the rest of the paper.

\begin{figure}
	\centering\includegraphics[width=0.45\textwidth] {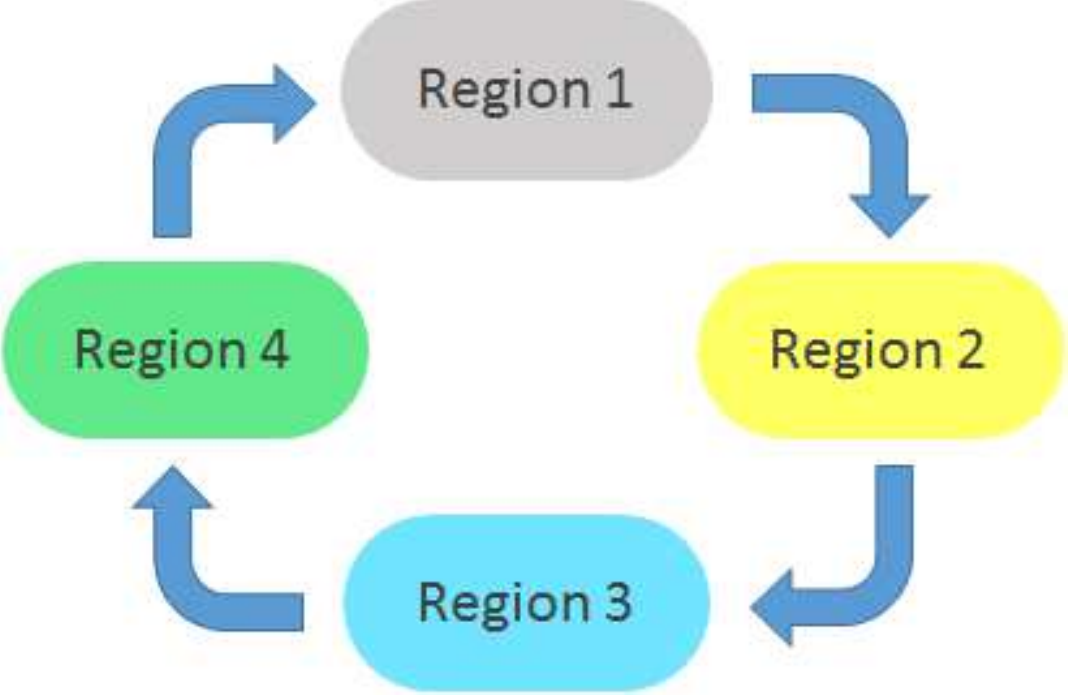}
	\caption{The cyclic evolution of solution trajectories of the system \eqref{eq:ghane}.}
	\label{fig:Cyclic_evolution}
\end{figure}

From the equations of motion~\eqref{eq:ghane} we see that the $x_3$-axis is also invariant under the flow. 
From the third component of \eqref{eq:ghane} we see that the plane $x_3=h $ is transversal to the flow at all points not contained in the intersection with the ellipsoid $a{x_1}^2 + b{x_2}^2 + c{x_3}^2 = r^2$.
The  ellipsoid $a{x_1}^2 + b{x_2}^2 + c{x_3}^2 = r^2$ is however not even away from its intersection with the plane $x_3=h $  transversal to the flow.  
Still we can use the structure of the NEValues summarized in  Table~\ref{tab:Qualitative} to infer that 
the system trajectories with initial conditions $x_3>0$ and not contained in the invariant $x_3$-axis visit the different region in the cyclic pattern illustrated in Figure~\ref{fig:Cyclic_evolution}. 
To this end first note that the NEValues show that trajectories cannot get permanently trapped in either of the regions 1, 2, 3 or 4.  
Let us consider trajectories going through the intersection of the nullclines $\dot{\rho}=0$ and $\dot{x}_3=0$ which is given by the ellipse $a{x_1}^2 + b{x_2}^2 + c{h}^2 = r^2$ in the plane $x_3=h$. The ellipse is illustrated in  Fig.~\ref{fig:Qualitative_analysis}c for the case $a>b$ (for $a<b$ relabel $x_1$ and $x_2$). Trajectories crossing the ellipse in the first or third quadrant evolve from region~3 to region~2. Trajectories crossing the ellipse in the second or fourth quadrant evolve from region~1 to region~4. At the four points on the ellipse contained on the $x_1$-axis or $x_2$-axis the vector field is tangent to the ellipse. At the two points on the ellipse located on the $x_1$-axis the trajectories evolve from region~3 to region~4. 
At the two points on the ellipse located on the $x_2$-axis the trajectories evolve from region~1 to region~2. 

Using the above results on how the nullclines $\dot{x}_3=0$ and  $\dot{\rho}=0$ are crossed by trajectories we can conclude the following on the evolution between the different regions.
All trajectories with initial conditions in region~4 will evolve directly (i.e. without visiting any other region in between) to region~1. 
Similarly, all trajectories with initial conditions in region~2 will evolve directly to region~3.
Trajectories with initial conditions in region~1 evolve directly to region~2, or move to region~2 after a finite number of visits to region~4. The latter option follows from the non-transversality of the boundary between regions~1 and 4 and the fact that $\rho$ is exponentially inreasing in regions~1 and 4.
Similarly, trajectories with  initial conditions in region~3 evolve directly to region~4, or move to region~4 after a finite number of visits to region~3. The latter option follows from the non-transversality of the boundary between regions~2 and 3 and the fact that $\rho$ is exponentially decreasing in regions~2 and 3.

By the arrangement of regions 1--4, it follows that the stretching of trajectories occurs along the  $x_3$-axis followed by the folding action and the cyclic pattern occurs repeatedly. On the other hand, the results of the qualitative analysis of the system trajectories summarized in Table~\ref{tab:Qualitative} and Figures~\ref{fig:Qualitative_analysis} and \ref{fig:Cyclic_evolution} suggest that the system will be globally bounded, which will be rigorously proved in the next subsection.

In summary, based on this qualitative analysis we expect that the synthesized system \eqref{eq:ghane} will exhibit chaotic behaviour for a suitable range of parameter values. Note that the proposed approach is essentially qualitative without any quantitative rigorous proof.
In the next section, we present a numerical analysis of the system which indeed suggests the occurrence of chaotic behaviour and we discuss some interesting features of the system.

\begin{table}
	\centering
	\caption{Qualitative behaviour of the nonlinear system \eqref{eq:ghane} based on an analysis of the NEValues.}
	\label{tab:Qualitative}
	\begin{tabular}{|l|l|l|ll|}
\hline
  & Region & Sign of ${\mathop{\rm Re}\nolimits} \left\{ {{\lambda }({\bf{x}})} \right\}$ & \multicolumn{2}{l|}{Qualitative behaviour}                                                               \\ \hline
 \multicolumn{1}{|c|}{\multirow{2}{*}{1}} & \multirow{2}{*}{$\begin{array}{l}a{x_1}^2 + b{x_2}^2 + c{x_3}^2 > {r^2}\; \\ {x_3}^2 > {h}^2\end{array}$} & ${\mathop{\rm Re}\nolimits} \left\{ {{\lambda _{1,2}}({\bf{x}})} \right\} >0$  & \multicolumn{1}{l|}{$\frac{d}{{dt}}\rho>0 $} & Increasing spiral      \\ \cline{3-5} 
		\multicolumn{1}{|c|}{}                          &                                                                                                                 & ${\mathop{\rm Re}\nolimits} \left\{ {{\lambda _3}({\bf{x}})} \right\} < 0$     & \multicolumn{1}{l|}{$\frac{d}{{dt}}(\left| {{x_3}} \right|) < 0$}             & Decreasing exponential \\ \hline
		\multirow{2}{*}{2}                      & \multirow{2}{*}{$\begin{array}{l}a{x_1}^2 + b{x_2}^2 + c{x_3}^2 > {r^2}\; \\ {x_3}^2 < {h}^2\end{array}$} & ${\mathop{\rm Re}\nolimits} \left\{ {{\lambda _{1,2}}({\bf{x}})} \right\} < 0$ & \multicolumn{1}{l|}{$\frac{d}{{dt}}\rho<0 $} & Decreasing spiral      \\ \cline{3-5} 
		&                                                                                                                 & ${\mathop{\rm Re}\nolimits} \left\{ {{\lambda _3}({\bf{x}})} \right\} < 0$     & \multicolumn{1}{l|}{$\frac{d}{{dt}}(\left| {{x_3}} \right|) < 0$}             & Decreasing exponential \\ \hline
		\multirow{2}{*}{3}                      & \multirow{2}{*}{$\begin{array}{l}a{x_1}^2 + b{x_2}^2 + c{x_3}^2 < {r^2}\; \\ {x_3}^2 < {h}^2\end{array}$} & ${\mathop{\rm Re}\nolimits} \left\{ {{\lambda _{1,2}}({\bf{x}})} \right\} < 0$ & \multicolumn{1}{l|}{$\frac{d}{{dt}}\rho<0 $} & Decreasing spiral      \\ \cline{3-5} 
		&                                                                                                                 & ${\mathop{\rm Re}\nolimits} \left\{ {{\lambda _3}({\bf{x}})} \right\} > 0$     & \multicolumn{1}{l|}{$\frac{d}{{dt}}(\left| {{x_3}} \right|) > 0$}             & Increasing exponential \\ \hline
		\multirow{2}{*}{4}                      & \multirow{2}{*}{$\begin{array}{l}a{x_1}^2 + b{x_2}^2 + c{x_3}^2 < {r^2}\; \\ {x_3}^2 >{h}^2\end{array}$}  & ${\mathop{\rm Re}\nolimits} \left\{ {{\lambda _{1,2}}({\bf{x}})} \right\} > 0$ & \multicolumn{1}{l|}{$\frac{d}{{dt}}\rho>0 $} & Increasing spiral      \\ \cline{3-5} 
		&                                                                                                                 & ${\mathop{\rm Re}\nolimits} \left\{ {{\lambda _3}({\bf{x}})} \right\} > 0$     & \multicolumn{1}{l|}{$\frac{d}{{dt}}(\left| {{x_3}} \right|) > 0$}             & Increasing exponential \\ \hline
	\end{tabular}
\end{table}

\subsection{Boundedness of system trajectories}

We note that  all orbits in the plane $x_3=0$ are attracted to the origin (see subsection~\ref{sec:equilibria}) which guarantees boundedness of 
trajectories with initial conditions in the plane $x_3=0$ in the forward time direction. 


Given the possible transport scenarios between the regions 1-4 discussed in the previous subsection the boundedness of system trajectories in general 
follows from the following proposition.

\begin{proposition}\label{prop:boundedness_general}
Trajectories with initial conditions in region~1 can enter region~2 only with a finite maximal value of $\rho$ (where this maximal value depends on the initial condition). 
\end{proposition}

This gives the boundedness of system trajectories in general because (as before we consider because of symmetry only the half $x_3>0$):
\begin{enumerate}
\item trajectories on the invariant $x_3$-axis have $\dot{x}_3<0$ for $x_3>r\sqrt{c}$,
\item regions 3 and 4 are bounded,
\item region 2 is bounded in the vertical direction from below and from above and in region~2 we have $\dot{\rho}<0$, and
\item region 1 can can only be entered from the bounded region $4$. 
\end{enumerate}

Let us now prove proposition~\ref{prop:boundedness_general}.

\begin{proof} (Proposition~\ref{prop:boundedness_general})
	We show that trajectories with initial conditions $(x_1(0),x_2(0),x_3(0))=(x_{1\,0},x_{2\,0},x_{3\,0})$ in region~1 will reach the plane $x_3=h$ in the forward time direction with a finite value $\rho_0$. 
	Let $x_{3\,0}=x_3(0)$ and $d=\text{min }\{a,b\}>0$. Then there exists a $\rho_1> \text{max }\{r/\sqrt{d},\sqrt{x_{1\,0}^2+x_{2\,0}^2}\}$ such that for all points $(x_1,x_2,x_3)$ with $h\le x_3\le x_{3\,0}$ and $\sqrt{x_1^2+x_2^2}=\rho\ge \rho_1$
\begin{equation}
	\begin{split}
	\frac{\dot{x}_3}{\dot{\rho}} = \frac{(r^2 -a x_1^2 -b x_2^2 - c x_3^2 )x_3}{(x_3^2-h^2)\rho} &\le  \frac{(r^2 -d \rho^2  )x_3}{(x_3^2-h^2)\rho} \\
	&\le  \frac{(r^2 -d \rho^2  )x_3}{(x_{3\,0}^2-h^2)\rho} \\ \label{eq:to_min_inf_for_rho_to_inf_general}
	&\le  \frac{(r^2 -d \rho^2  )h}{(x_{3\,0}^2-h^2)\rho}\\ 
	&<-1.
	\end{split}
\end{equation}	
	The existence of $\rho_1$ follows from the last but one expression in \eqref{eq:to_min_inf_for_rho_to_inf_general} going to $-\infty$ as $\rho\to\infty$.
	As $\dot{x}_3<0$ in region~1 the vertical variation of the tajectory in the forward time direction is equal to  $x_{3\,0}-h$  before the orbit reaches the plane $x_3=h$. 
	The orbit can then depart in the forward time direction no further from the $x_3$-axis than
	$\rho_0=  \rho_1 + x_{3\,0}-h$ before it reaches the plane $x_3=h$.  	
\end{proof}


\section{Dynamical analysis of the candidate system}

In this section we study the dynamics of the system \eqref{eq:ghane}. 
We start by studying the bifurcations of equilibria and periodic orbits. Numerical simulations suggest that chaotic attractors appear after a cascade of period doubling bifurcations. These chaotic attractors are presumably of H\'enon-like type which means that they are the closure of the unstable manifold of a saddle periodic orbit.


\subsection{Equilibria and their stability}
\label{sec:equilibria}

For $x_3 \geq 0$, the system has the following equilibria:
\[
	O = (0,0,0), \quad
	Z = (0,0,r/\sqrt{c}).
\]
The eigenvalues of the Jacobi matrix $J=D \bf{f}$ evaluated at $O$ are given by
\[
	\lambda_{O,1} = -h^2 + \omega j, \quad
	\lambda_{O,2} = -h^2 - \omega j, \quad
	\lambda_{O,3} = r^2.
\]
Under the assumption that $h, r \neq 0$ the stable and unstable manifolds of $O$ are given by
\[
	\begin{split}
	W^s(O) & = \{(x_1, x_2, 0) \in \mathbb{R}^3\}, \\
	W^u(O) & = \{(0, 0, x_3) \in \mathbb{R}^3 \,:\, 0 < |x_3| < r/\sqrt{c}\}.
	\end{split}
\]
The eigenvalues of the matrix $J$ evaluated at $Z$ are given by:
\[
	\lambda_{Z,1} = r^2/c-h^2 + \omega j, \quad
	\lambda_{Z,2} = r^2/c-h^2 - \omega j, \quad
	\lambda_{Z,3} = -2r^2.
\]
Note that the complex eigenvalue pair $(\lambda_{Z,1}, \lambda_{Z,2})$ crosses the imaginary axis when $h = r / \sqrt{c}$. This implies that, under suitable non-degeneracy conditions, the equilibrium $Z$ becomes unstable through a Hopf bifurcation which gives birth to a stable periodic orbit (see Kuznetsov (2004)).

For $a=b$ and $r^2/c - h^2>0$ it is straightforward to verify that
\begin{equation}
\label{eq:periodic-orbits}
	\begin{split}
	x_1(t) & = \rho_0\cos(\omega t), \\
	x_2(t) & = \rho_0\sin(\omega t), \\
	x_3(t) & = h,
	\end{split}
\end{equation}
where $\rho_0 = \sqrt{(r^2-ch^2)/a}$, is a periodic orbit of the system \eqref{eq:ghane}. Note that for $h = r / \sqrt{c}$ this orbit coalesces with the equilibrium $Z$. This suggests that for $a=b$ the periodic orbit in equation \eqref{eq:periodic-orbits} indeed arises through a Hopf bifurcation of the equilibrium $Z$.

\begin{figure}
	\centering
	\includegraphics[width=0.48\textwidth]{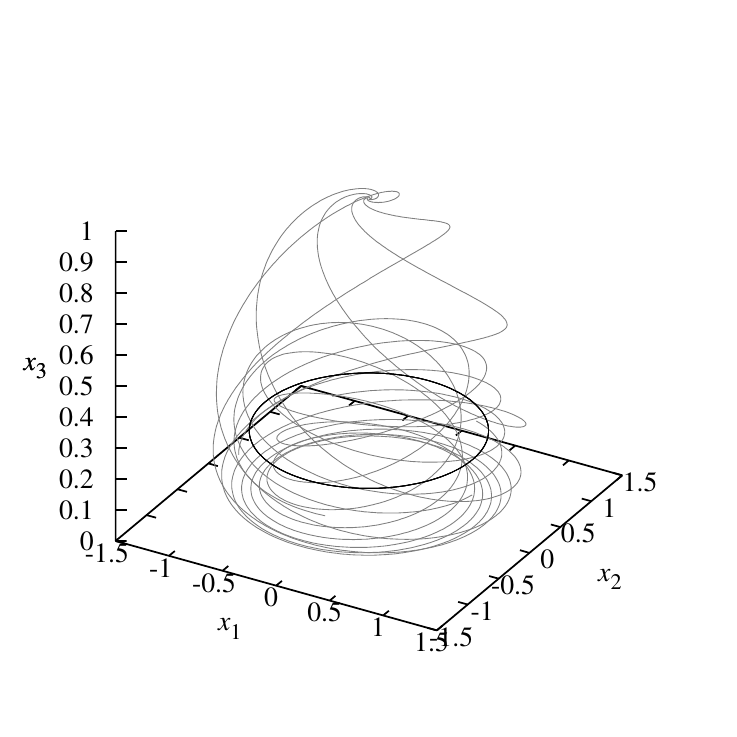}
	\includegraphics[width=0.48\textwidth]{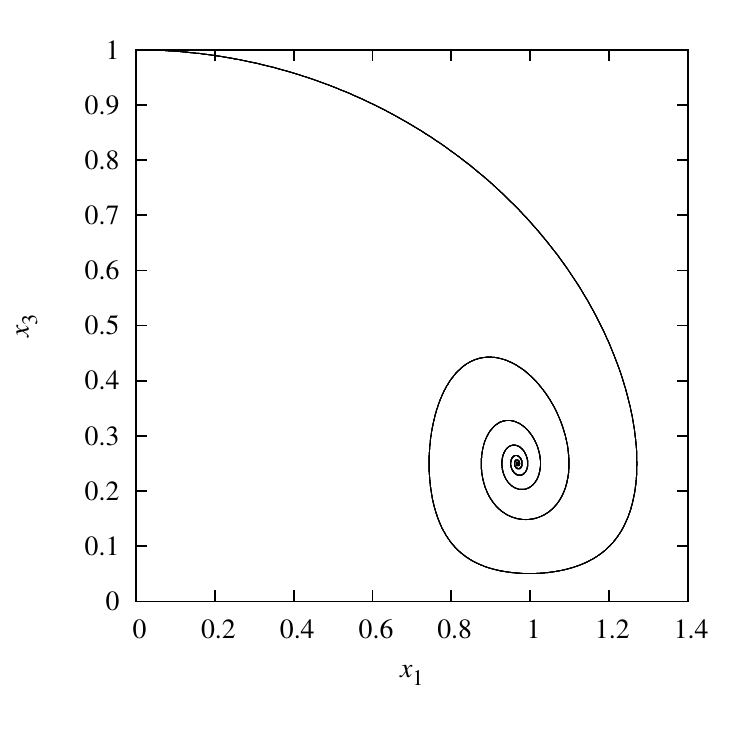}
	\caption{Left: Five orbits on the unstable manifold of the equilibrium $Z$ for the parameters $(a,b,c,h,r,\omega,c)=(1,1,1,0.25,1,1)$ (grey) and a stable periodic orbit (black). Right: cross section of the unstable manifold of the equilibrium $Z$.}
	\label{fig:manifold}
\end{figure}

Under the assumptions that $r \neq 0$ and $r^2/c - h^2>0$ it follows that the stable manifold of $Z$ is given by
\[
	W^s(Z) = \{(0, 0, x_3) \in \mathbb{R}^3 \,:\, x_3 > 0\}.
\]
The 2-dimensional unstable manifold of $Z$ cannot be computed analytically, but the linearization of the system \eqref{eq:ghane} at $Z$ shows that the unstable manifold is tangent to the plane $\{(x_1,x_2,r/\sqrt{c}) \,:\, x_1, x_2 \in \mathbb{R}\}$. Figure~\ref{fig:manifold} shows a numerical approximation of $W^u(Z)$, which suggests that this manifold is part of the stable manifold of a periodic orbit.

\subsection{Periodic orbits and their bifurcations}

The periodic solutions of \eqref{eq:ghane} can be studied in terms of a so-called Poincar\'e return map (see Guckenheimer \& Holmes (1983)). The idea is to study the intersections of orbits of \eqref{eq:ghane} with a plane that is transversal to the vector field. Consider the following set:
\[
	\Sigma = \{(x_1,0,x_3) \in \mathbb{R}^3 \,:\, x_1>0, \, x_3>0\}.
\]
We define the Poincar\'e map $P : \Sigma \to \Sigma$ as follows. If $(x_1,0,x_3) \in \Sigma$, then $P(x_1,0,x_3)$ is defined by integrating equation \eqref{eq:ghane} for $2\pi/\omega$ units of time. From equation \eqref{eq:ghane-cylinder} it follows  that indeed $P(\Sigma)\subset \Sigma$. In addition, the existence and uniqueness theorems for differential equations imply that the map $P$ is a diffeomorphism. A point ${\bf x} \in \Sigma$ is called a period-$n$ point of $P$ if $P^n({\bf x}) = {\bf x}$. Such points correspond to periodic orbits of \eqref{eq:ghane} which make $n$ turns around the $x_3$-axis. Period-1 points of $P$ are also referred to as fixed points of $P$.

\begin{figure}
	\centering
	\includegraphics[width=0.49\textwidth]{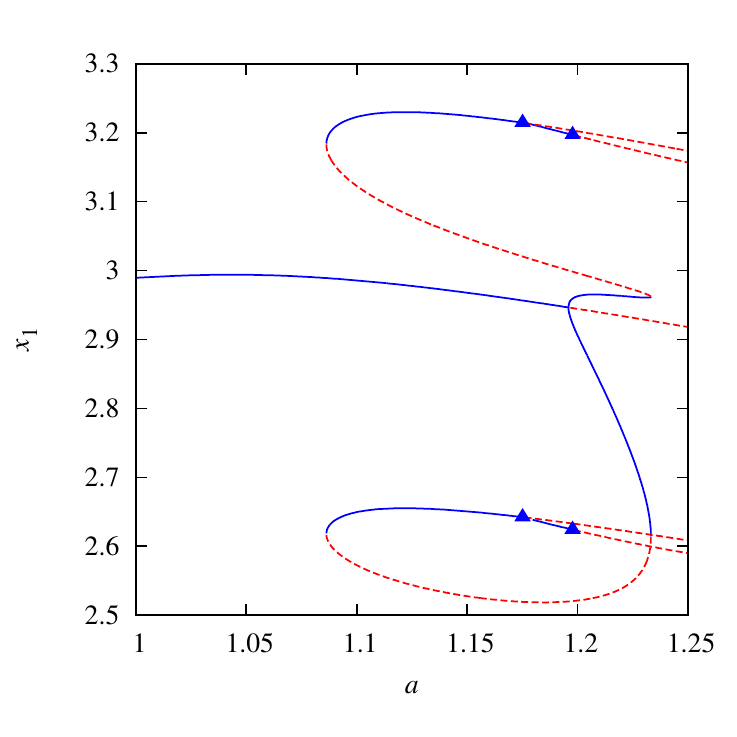}
	\caption{Bifurcation diagram of fixed points of the Poincar\'e map $P : \Sigma \to \Sigma$ as a function of the parameter $a$. The other parameters are fixed at $(b,c,h,r,\omega)=(1,1,0.25,3,1)$. Blue solid lines indicate stable branches and red dashed lines represent unstable branches. Triangles denote period doubling bifurcations. Note that multiple stable fixed points can coexist for the same parameter values.}
	\label{fig:bifurcations-periodic}
\end{figure}

For $a=b$ and $r^2/c - h^2>0$, the point $(\rho_0, 0, h)$, with $\rho_0 = \sqrt{(r^2-ch^2)/a}$, is a fixed point of the map $P$. This fixed point corresponds to the periodic solution given in 
equation~\eqref{eq:periodic-orbits}. Using the numerical continuation software package AUTO-07P (see Doedel \& Oldeman (2007)), we have computed the bifurcation diagram for this fixed point shown in Figure~\ref{fig:bifurcations-periodic}. The parameter $a$ is used as the continuation parameter; the other parameters are fixed at $(b,c,h,r,\omega)=(1,1,0.25,3,1)$. The fixed point is stable up to $a \approx 1.196$ where it loses stability in a supercritical pitchfork bifurcation. From the pitchfork bifurcation two stable fixed points emanate which lose (resp.\ regain) stability at saddle-node bifurcations for $a \approx 1.233$ (resp.\ $a\approx 1.086$). After that the two branches undergo a period doubling bifurcation at $a \approx 1.175$. This leads to the coexistence of two stable period-2 points.

For $a\approx 1.197$ the stable period-2 points lose stability through a period doubling bifurcation which leads to the coexistence of two stable period-4 points. This suggests that an infinite cascade of period doubling bifurcations occurs when $a$ increases. In principle, the next period doubling bifurcations can be obtained by means of numerical continuation. However, in a period doubling cascade the distances between successive period doublings asymptotically scale with the Feigenbaum constant $\delta \approx 4.669$ (see Guckenheimer \& Holmes (1983)). This implies that prohibitively small step sizes are needed to detect the bifurcations by means of continuation. Bifurcations can go undetected when the step size is too large.

A more practical way of obtaining an overview of the dynamics of the Poincar\'e map $P$ is to use brute force iteration. We increase the value of $a$ from $1.197$ up to $1.205$ in $1000$ steps. For each value of $a$ we compute $600$ iterates of $P$ and plot the $x_1$-coordinates of the last $100$ computed points as a function of $a$. The final point of the last attractor serves as an initial condition for the next loop. The starting points are the two stable period-2 points $(2.633, 0.00129)$ and $(3.203, 0.03657)$. The bifurcation diagrams for these points are shown in Figure \ref{fig:bifurcations-periodic}. This figure suggests that indeed each of the two points bifurcates through an infinite cascade of period doublings. In turn this leads to the coexistence of two chaotic attractors of which the structure will be discussed in the next section.

The coexistence of two or more attractors in a dynamical systems is referred to as \emph{multi-stability}. This phenomenon often arises due to symmetries of the system (see Lai \& Chen (2016)) and in particular due to the presence of pitchfork bifurcations (see Van Kekem \& Sterk (2017) and (2018b)) as is the case in the present paper. A different mechanism by which multi stability can occur is due to the presence of codimension-2 bifurcations, such as double-Hopf bifurcations (see Van Kekem \& Sterk (2018a)). For an overview of the wide range of applications of multi-stability in different disciplines of science, see Feudel (2008).

\begin{figure}
	\centering
	\includegraphics[width=0.49\textwidth]{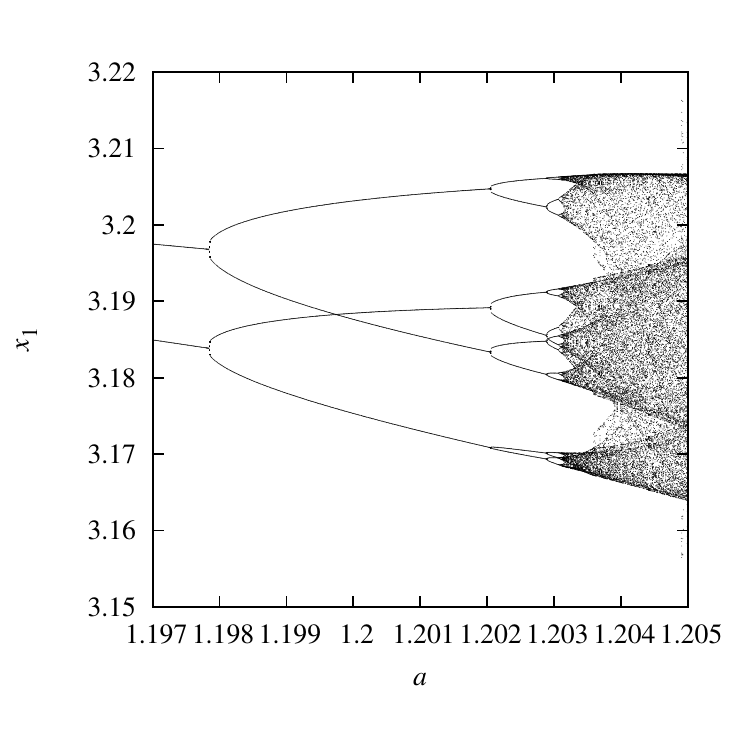}
	\includegraphics[width=0.49\textwidth]{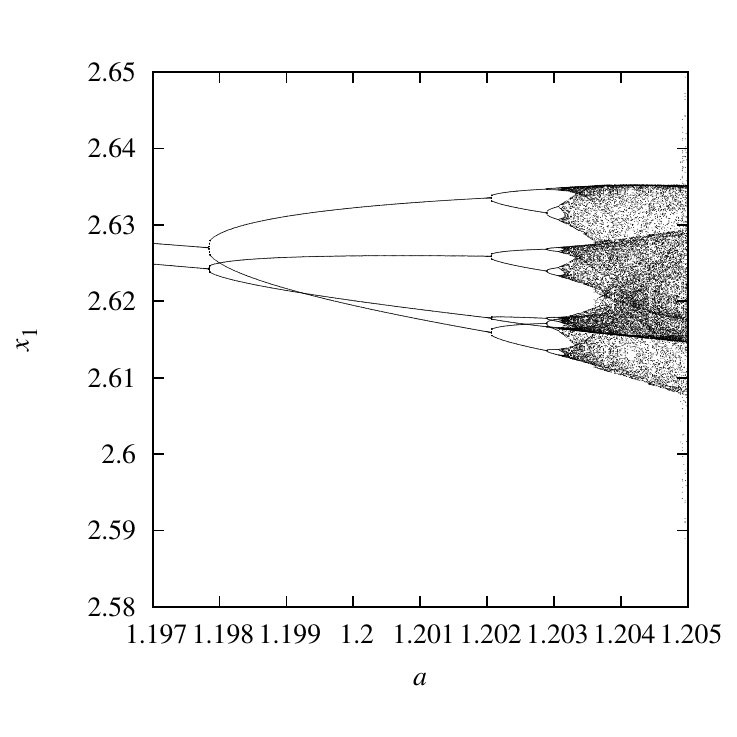}
	\caption{Bifurcation diagram of two stable period-2 points of the Poincar\'e map $P : \Sigma \to \Sigma$. For each value of $a$ the map $P$ is iterated 500 times and the last 100 iterates are plotted. The parameters $(b,c,h,r,\omega)=(1,1,0.25,3,1)$ are fixed.}
	\label{fig:cigars}
\end{figure}


\subsection{Chaotic dynamics}

\begin{figure}
 	\centering
 	\includegraphics[width=0.49\textwidth]{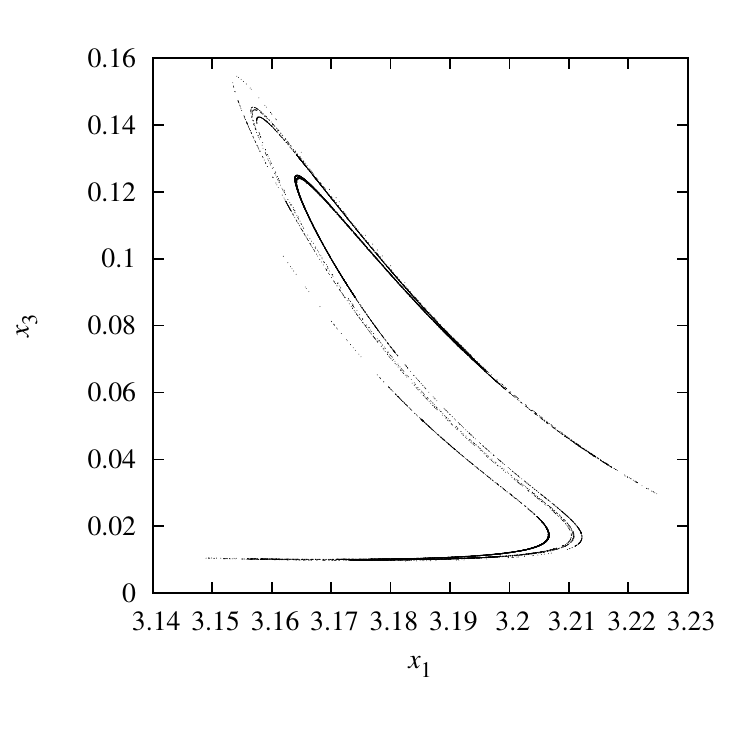}
 	\includegraphics[width=0.49\textwidth]{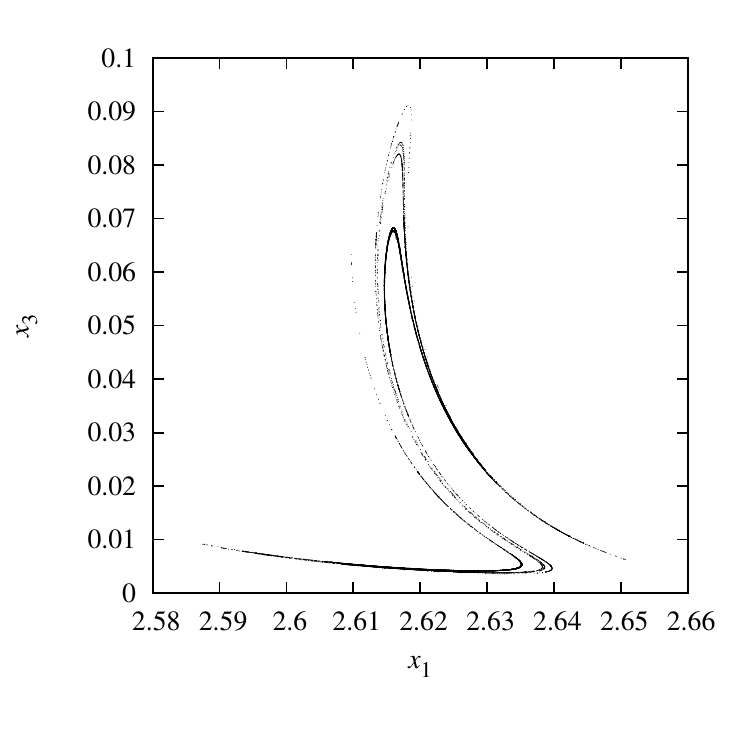}
 	\caption{Two coexisting chaotic attractors of the Poincar\'e map $P : \Sigma \to \Sigma$ for the parameter values $(a,b,c,h,r,\omega)=(1.205,1,1,0.25,3,1)$. The corresponding attractors for the system \eqref{eq:ghane} in $\mathbb{R}^3$ are shown in Figure \ref{fig:chaotic-attractors-flow}.}
 	\label{fig:chaotic-attractors}
 \end{figure}

Figure \ref{fig:chaotic-attractors} shows two chaotic attractors of the Poincar\'e map $P : \Sigma \to \Sigma$ detected after the period doubling cascade. Note that these attractors coexist for the same parameter values. The attractors have the appearance of a ``fattened curve'' which makes them qualitatively similar to the well-known attractor of the H\'enon map (see H{\'e}non (1976)). In fact, for the latter map it was proven by  Benedicks \& Carleson (1991) that for a set of parameter values with positive Lebesgue measure the attractor is the closure of the unstable manifold of a saddle fixed point.

\begin{figure}
 	\centering
 	\includegraphics[width=0.49\textwidth]{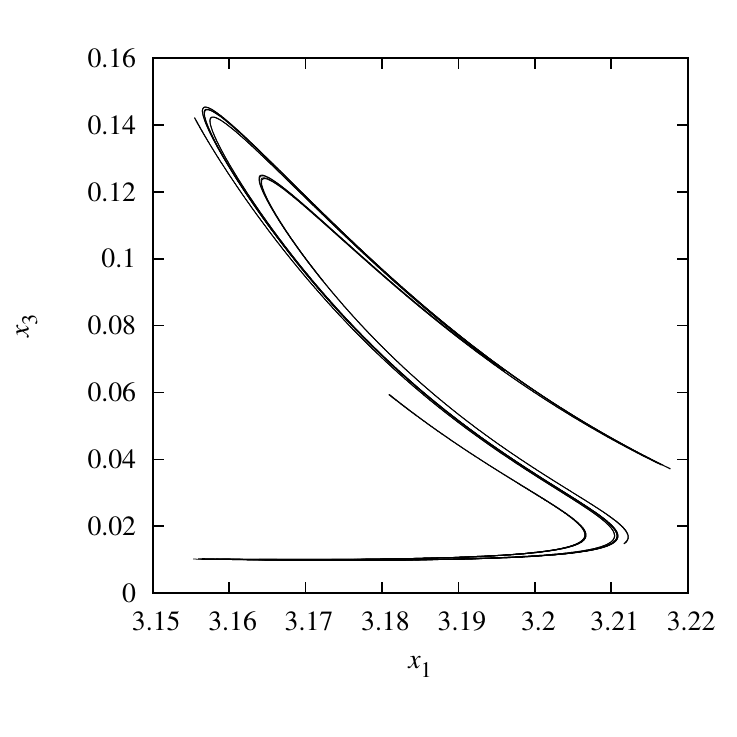}
 	\includegraphics[width=0.49\textwidth]{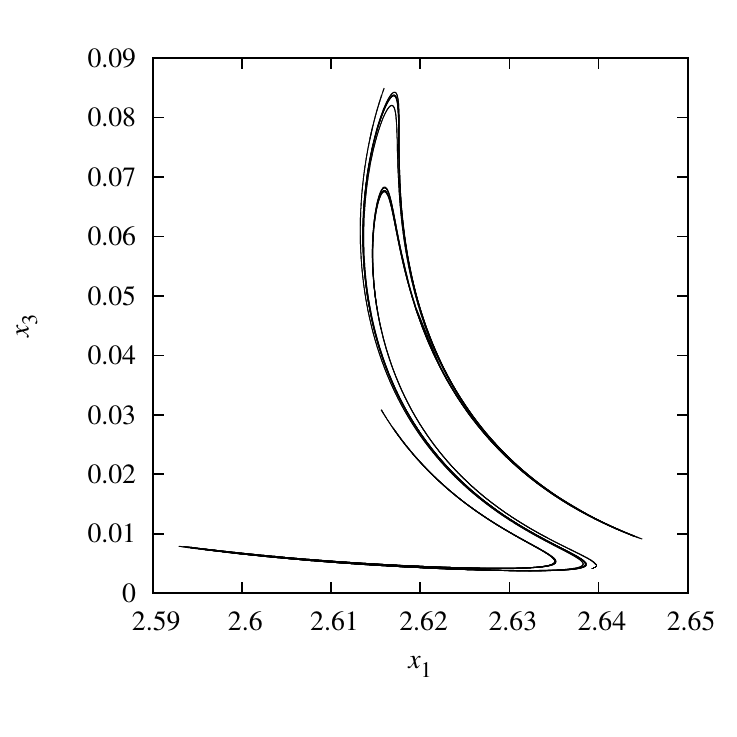}
 	\caption{Unstable manifold of two different saddle fixed points of the Poincar\'e map $P : \Sigma \to \Sigma$ for the parameter values $(a,b,c,h,r,\omega)=(1.205,1,1,0.25,3,1)$. Note the striking resemblance with the attractors shown in Figure \ref{fig:chaotic-attractors}.}
 	\label{fig:manifolds}
 \end{figure}

By numerical continuation we obtained two saddle fixed points of the Poincar\'e map $P : \Sigma \to \Sigma$ for the parameter values $(a,b,c,h,r,\omega)=(1.205,1,1,0.25,3,1)$. We computed the unstable manifolds of these fixed points by means of techniques based on iterating fundamental domains described in Broer \& Takens (2010) and Sim\'o (1990). Their unstable manifolds are shown in Figure \ref{fig:manifolds}. Note the striking resemblance with the attractors shown in Figure \ref{fig:chaotic-attractors}. We therefore conjecture that these attractors are in fact the closure of the manifolds shown in Figure~\ref{fig:manifolds}. This implies that the corresponding chaotic attractors for the system \eqref{eq:ghane} are the closure of the unstable manifold of a saddle periodic orbit.

 \begin{figure}
 	\centering
 	\includegraphics[width=0.49\textwidth]{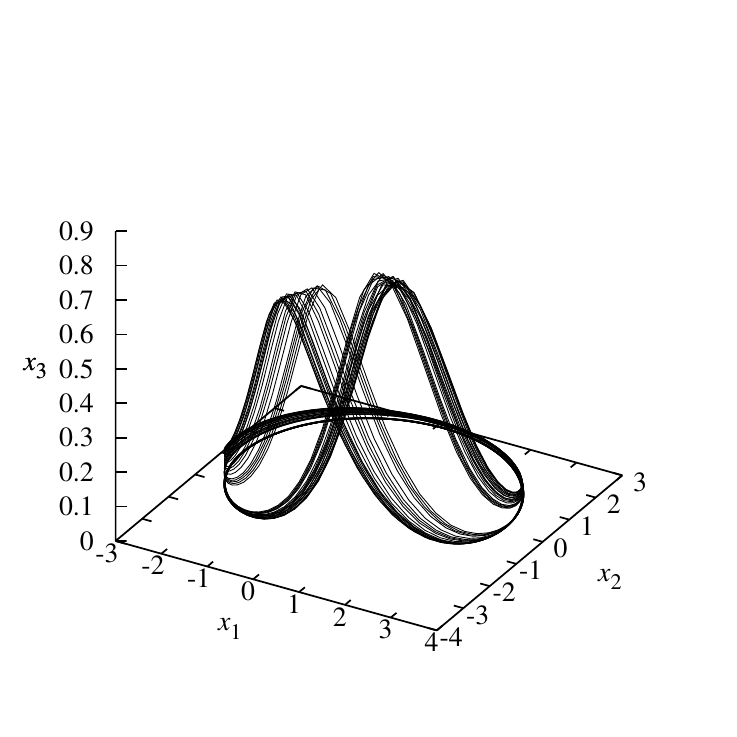}
 	\includegraphics[width=0.49\textwidth]{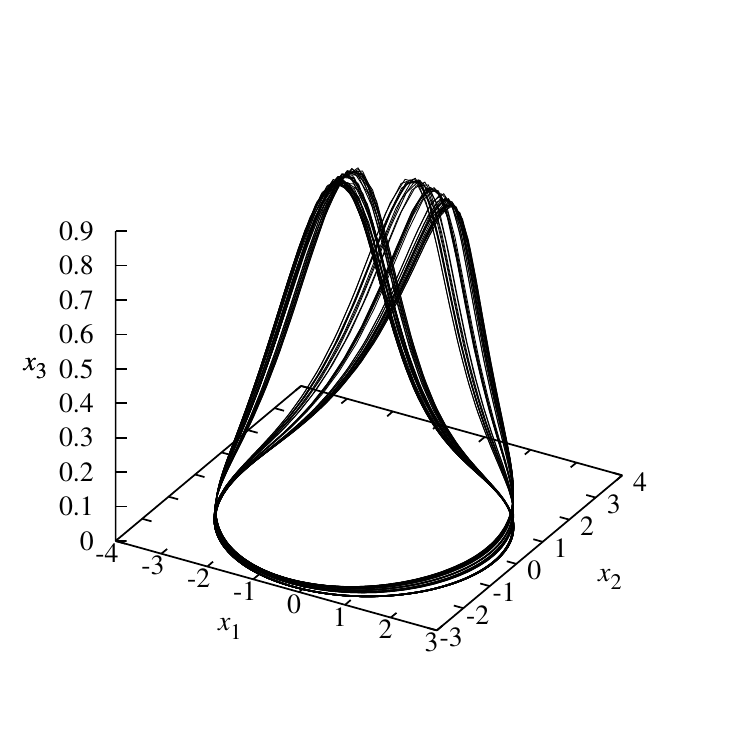}
 	\caption{Two coexisting chaotic attractors of the system \eqref{eq:ghane} for the parameter values $(a,b,c,h,r,\omega)=(1.205,1,1,0.25,3,1)$. Their corresponding Poincar\'e sections are shown in Figure \ref{fig:chaotic-attractors}.}
 	\label{fig:chaotic-attractors-flow}
 \end{figure}


\section{Control design} 

This section presents the possibility of using the eigen-structure analysis in nonlinear control design. First, a nonlinear state feedback controller is constructed to stabilize a chaotic system by the help of a NEValues assignment. Secondly, a synchronizing controller is obtained through a master-slave formalism.


\subsection{Chaos control}

The chaotic system \eqref{eq:ghane} can be controlled just by one control input. 
In fact we can control the system in such a way that the origin becomes a global attractor by a single input $u$ exerted on the $x_3$ component of   the system \eqref{eq:ghane} according to
\begin{equation}
	\label{eq:ghaneControl}
	\begin{split}
	\dot{x}_1 & = ({x_3}^2 - {h^2}){x_1} - \omega {x_2}, \\
	\dot{x}_2 & = \omega {x_1} + ({x_3}^2 - {h^2}){x_2}, \\
	\dot{x}_3 & = ({r^2} - a{x_1}^2 - b{x_2}^2 - c{x_3}^2){x_3} + u.
	\end{split}
\end{equation}
The control function $u : \mathbb{R} \to \mathbb{R}$ makes the origin asymptotically stable if the following condition is satisfied:
\begin{equation}
	\lambda_{3cl}({\bf{x}}) < 0 \quad\text{for all}\quad {\bf{x}} \in \mathbb{R}^3 \setminus \{{\bf{0}}\}.
\end{equation}

\begin{figure}
	\centering
	\begin{subfigure}[b]{0.49\textwidth}
		\includegraphics[width=\textwidth] {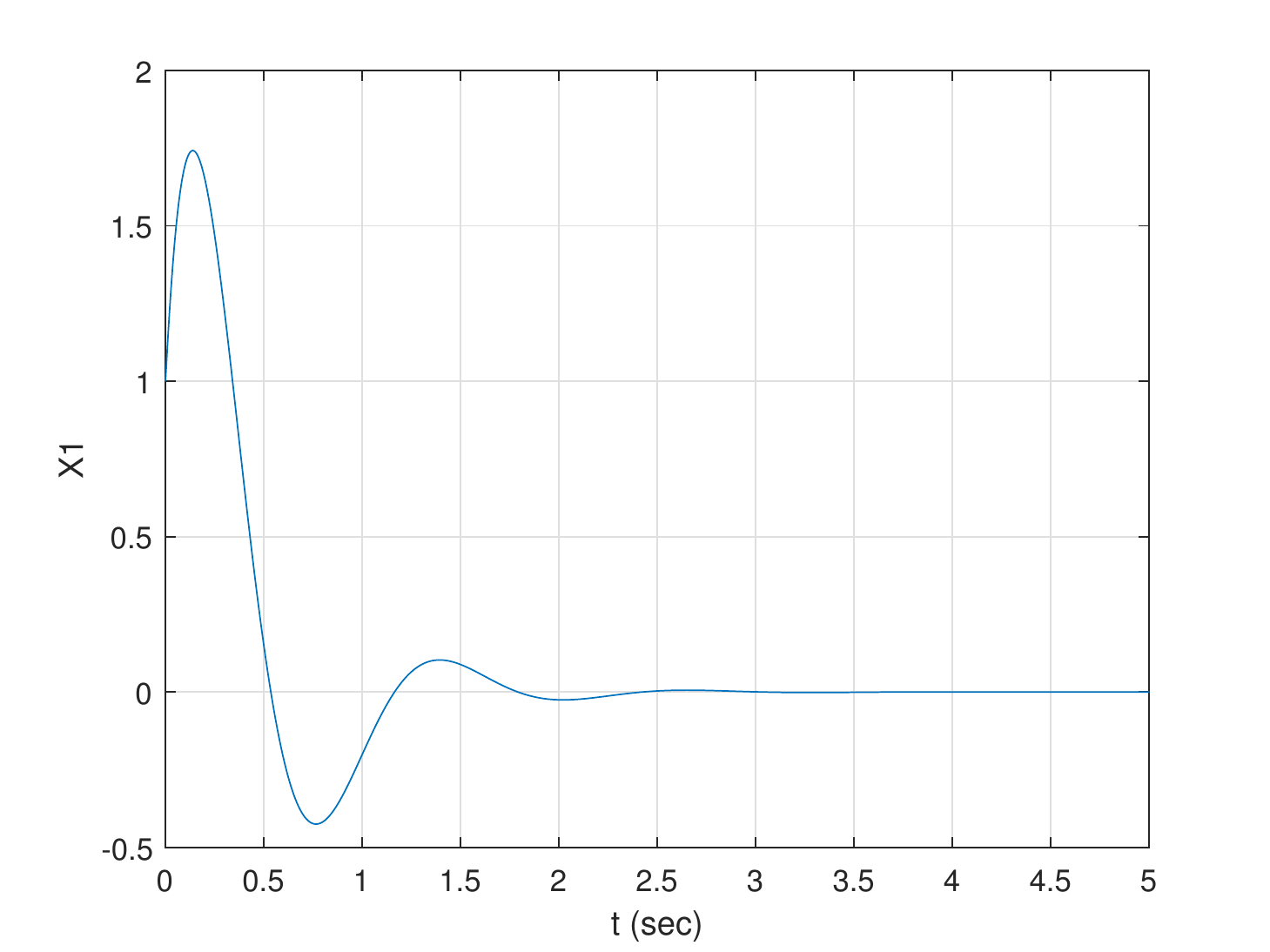}
			\end{subfigure}
	\begin{subfigure}[b]{0.49\textwidth}
		\includegraphics[width=\textwidth]{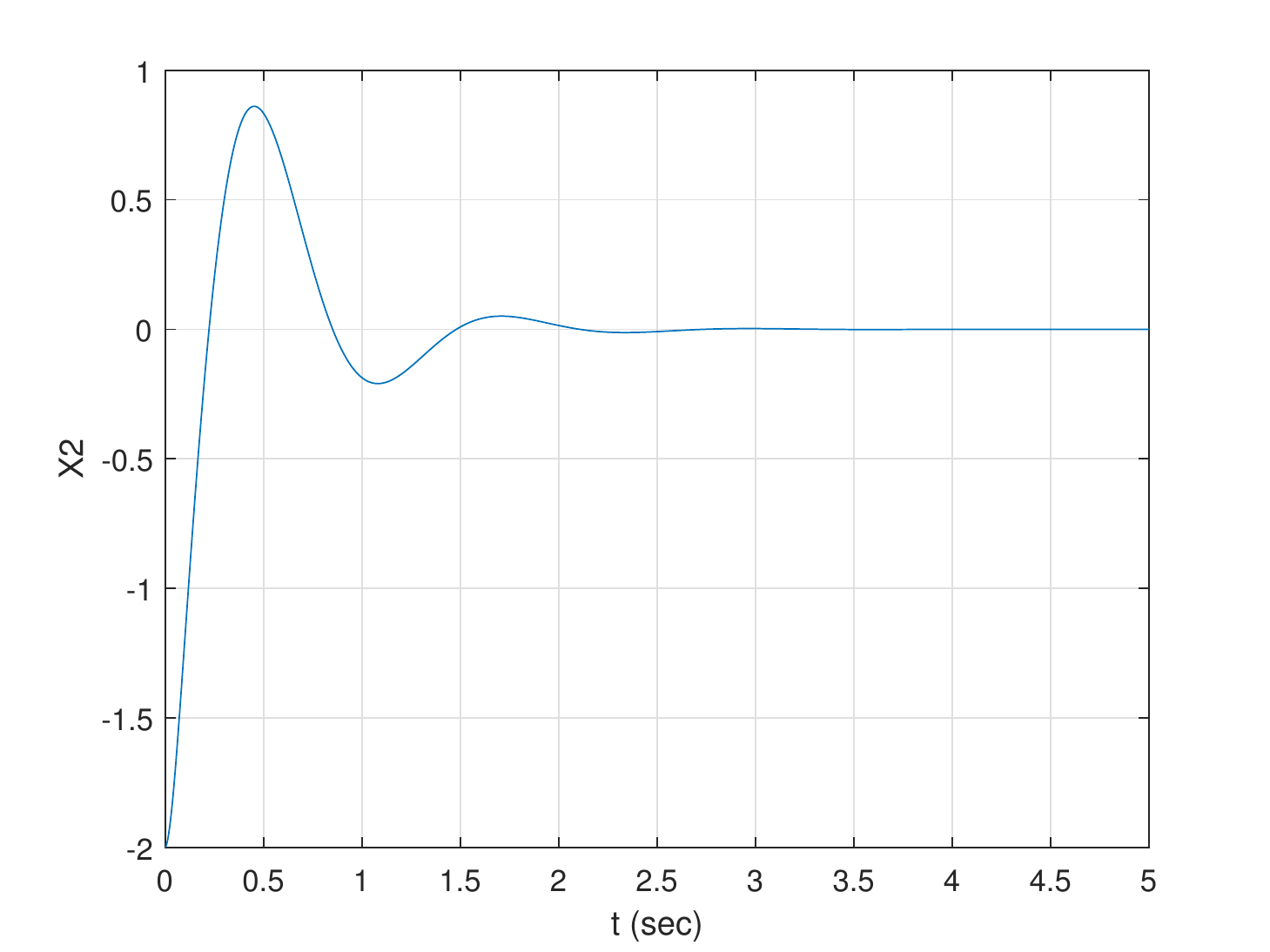}
	\end{subfigure}
	\begin{subfigure}[b]{0.49\textwidth}
		\includegraphics[width=\textwidth]{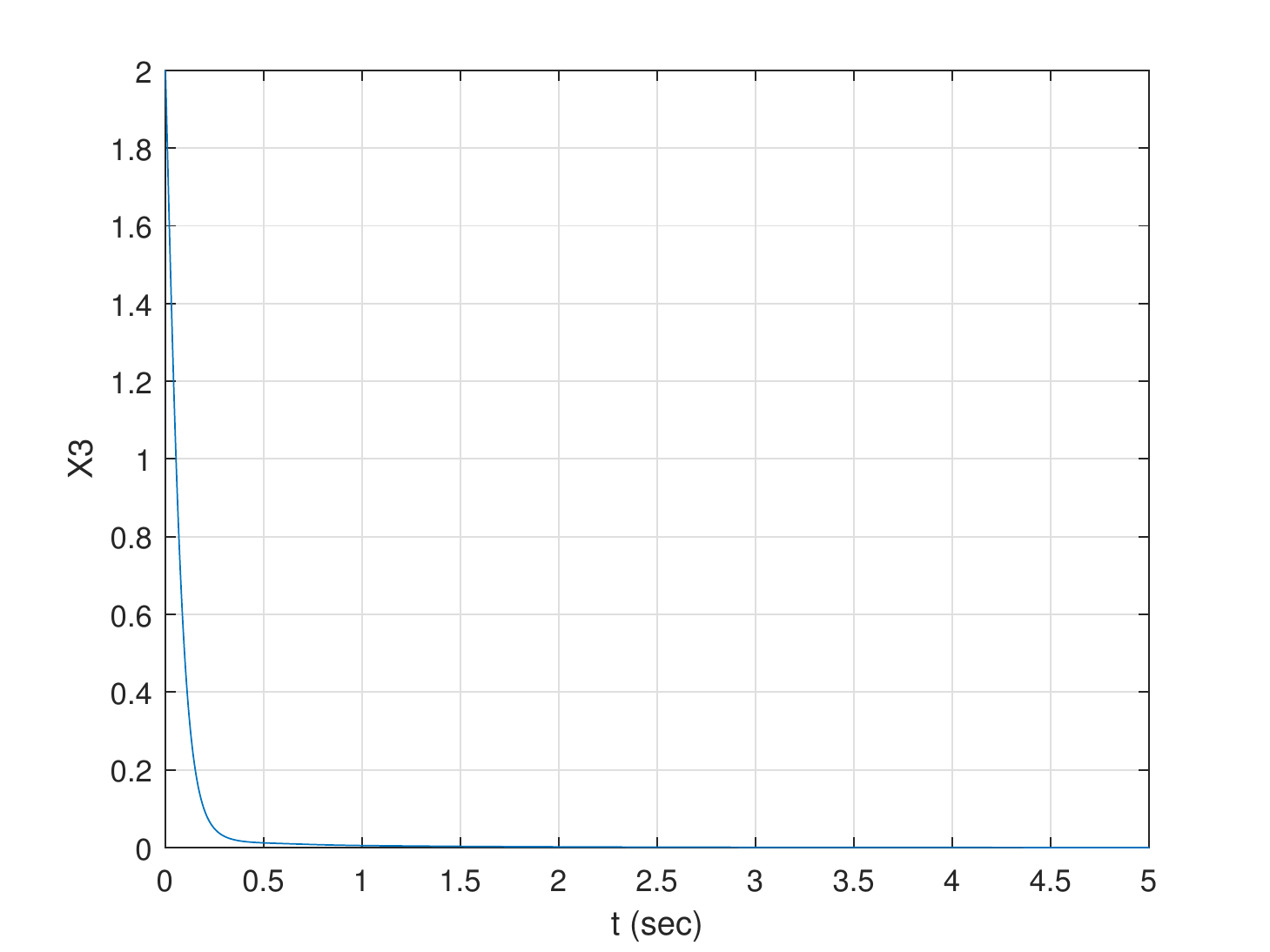}
	\end{subfigure}
	\begin{subfigure}[b]{0.49\textwidth}
		\includegraphics[width=\textwidth]{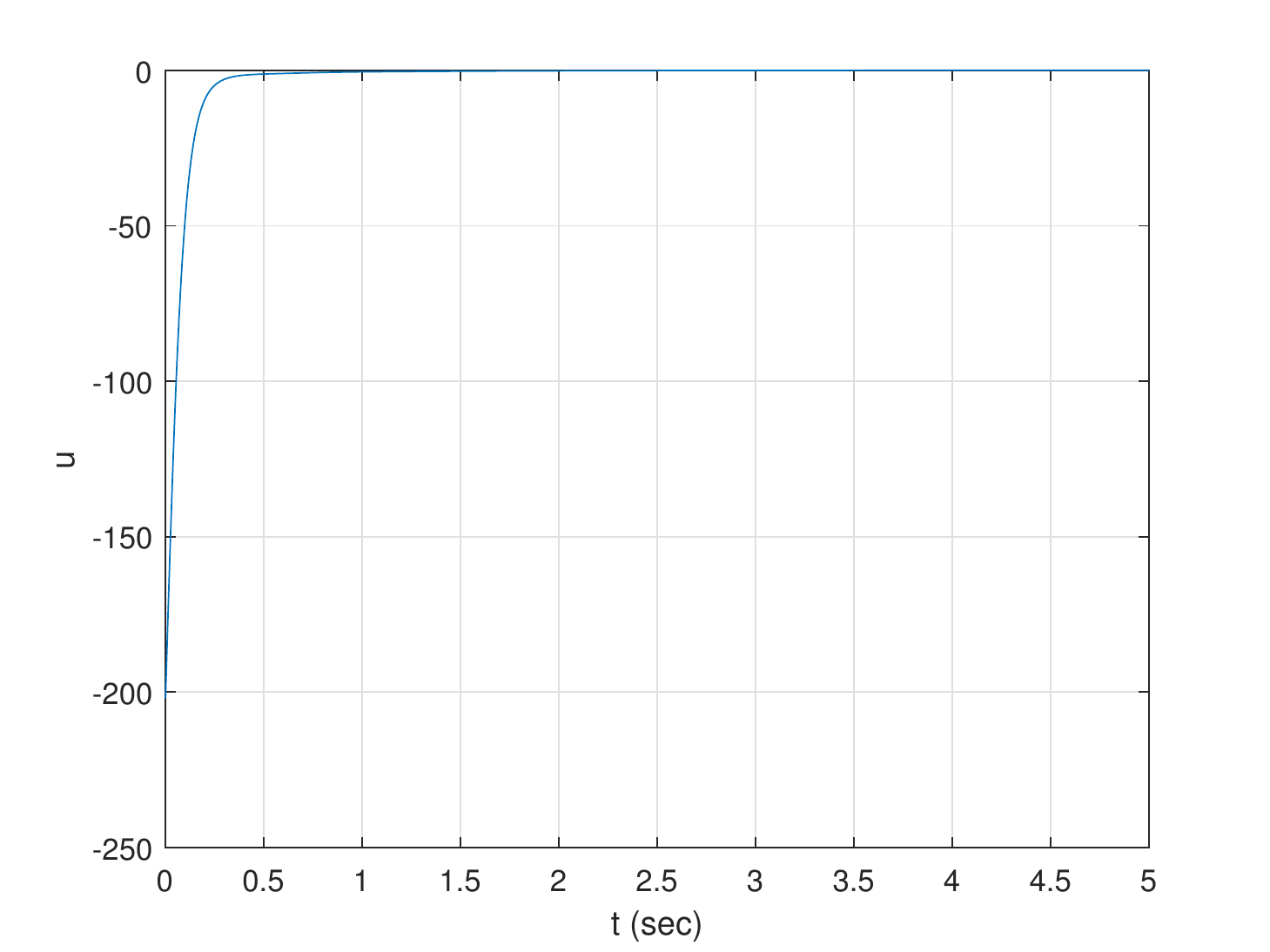}
	\end{subfigure}
	\caption{Closed loop results of system \eqref{eq:ghaneControl} with $u=-Kr^{2}x_{3}$ and $ K=1.1 $ for the parameter values $(a,b,c,h,r,\omega)=(5,1,0.1,1.5,10,5)$.}
	\label{fig:Control2}
\end{figure}

From the eigen-structure analysis of the system depicted in Table~1 we see that
even though the states $x_1$ and $x_2$ are not accessed by the input, the chaotic system can still be controlled by means of a simple state feedback of the form $u=-Kr^{2}x_{3}$ with $K>1$. The closed loop system obtained by applying this controller is
\[
	\begin{split}
	\dot{x}_1 & = ({x_3}^2 - {h^2}){x_1} - \omega {x_2}, \\
	\dot{x}_2 & = \omega {x_1} + ({x_3}^2 - {h^2}){x_2}, \\
	\dot{x}_3 & = ( (1-K)r^{2}- a{x_1}^2 - b{x_2}^2 - c{x_3}^2){x_3},
	\end{split}
\]
which is again in PL from with the old NEValues ${\lambda _{1,2}}_{cl}({\bf{x}})={\lambda _{1,2}}({\bf{x}})$ and the new NEValue  ${\lambda _3}_{cl}({\bf{x}}) =  (1-K)r^{2}- a{x_1}^2 - b{x_2}^2 - c{x_3}^2$. The simulation results of this controlled system are illustrated in Figure~\ref{fig:Control2}.

\subsection{Synchronization}

In this section we synchronize a pair of chaotic systems, which consist of a master system given by
\[
	\begin{split}
		{\dot x}_{m1} & = ({x_{m3}}^2 - {h^2}){x_{m1}} - \omega {x_{m2}}, \\
		{\dot x}_{m2} & = \omega {x_{m1}} + ({x_{m3}}^2 - {h^2}){x_{m2}}, \\
		{\dot x}_{m3} & = ({r^2} - a{x_{m1}}^2 - b{x_{m2}}^2 - c{x_{m3}}^2){x_{m3}}, 
	\end{split}
\]
and a slave system given by
\[
	\begin{split}
		{\dot x}_{s1} & = ({x_{s3}}^2 - {h^2}){x_{s1}} - \omega {x_{s2}}+u_{1}, \\
		{\dot x}_{s2} & = \omega {x_{s1}} + ({x_{s3}}^2 - {h^2}){x_{s2}}+u_{2}, \\
		{\dot x}_{s3} & = ({r^2} - a{x_{s1}}^2 - b{x_{s2}}^2 - c{x_{s3}}^2){x_{s3}} + u_{3}.
	\end{split}
\]
The equations for the error signal ${\bf{e}} = {{\bf{x}}_s} - {{\bf{x}}_m}$ are given by
\begin{equation}
	\label{eq:error_system}
	\begin{split}
		{\dot e}_1 & =  - {h^2}{e_1} - \omega {e_2} + {x_{s3}}^2{x_{s1}} - {x_{m3}}^2{x_{m1}}+u_{1}, \\
		{\dot e}_2 & = \omega {e_1} - {h^2}{e_2} + {x_{s3}}^2{x_{s2}} - {x_{m3}}^2{x_{m2}}+u_{2}, \\
		{\dot e}_3 & = {r^2}{e_3} - a({x_{s3}}{x_{s1}}^2 - {x_{m3}}{x_{m1}}^2) - b({x_{s3}}{x_{s2}}^2 - {x_{m3}}{x_{m2}}^2) - c({x_{s3}}^3 - {x_{m3}}^3) + u_{3}.
	\end{split}
\end{equation}

Unfortunately, the synchronization of the master and slave system cannot be achieved by only one control input. Instead, we need  three independent control inputs to guarantee that the origin of the system \eqref{eq:error_system} is asymptotically stable. Hence, we consider the controller
\[
\begin{split}
	{u_1} & = - {x_{s3}}^2{x_{s1}} + {x_{m3}}^2{x_{m1}}, \\
	{u_2} & = - {x_{s3}}^2{x_{s2}} + {x_{m3}}^2{x_{m2}}, \\
	{u_3} & = - {r^2}{e_3},
\end{split}
\]
which gives the following closed loop system for the error equations:
\[
\begin{split}
	{\dot e}_1 & = - {h^2}{e_1} - \omega {e_2}, \\
	{\dot e}_2 & = \omega {e_1} - {h^2}{e_2}, \\
	{\dot e}_3 & = - a({x_{s3}}{x_{s1}}^2 - {x_{m3}}{x_{m1}}^2) - b({x_{s3}}{x_{s2}}^2 - {x_{m3}}{x_{m2}}^2) - c({x_{s3}}^3 - {x_{m3}}^3).
\end{split}
\]

\begin{figure}
	\centering\includegraphics[width=\textwidth] {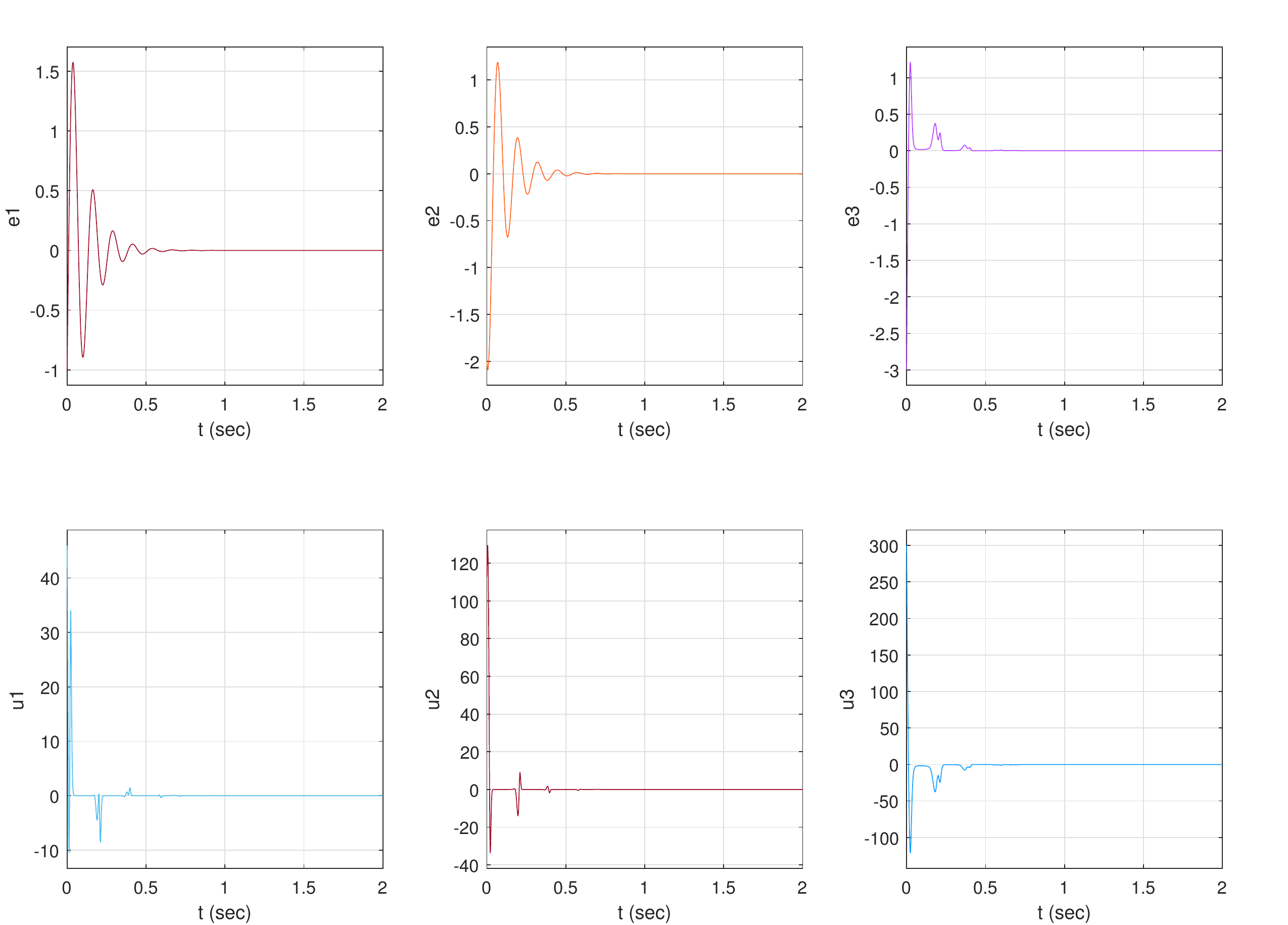}
	\caption{Synchronization simulation results for the parameter values $(a,b,c,h,r,\omega,c)=(5,1,0.1,4,10,50)$}
 	\label{fig:synchronization}	
\end{figure}

This means that regardless of ${e_3}$, the ${e_1} - {e_2}$ components approach zero in a spiralling manner  with the exponential rate of $ h^2 $. The remaining error dynamic in ${e_3}$ is then
\[
	{\dot e}_3 =  - \left( {a{x^{*}_1}^2 + b{x^{*}_2}^2 + c({x_{s3}}^2 + {x_{s3}}{x_{m3}} + {x_{m3}}^2)} \right){e_3}
\]
in which $x^{*}_1$ and $x^{*}_2$ denote the synchronized values of ${x_1}$ and ${x_2}$  assuming that ${e_{1,2}} $ has saturated at $0$. The dynamic of the ${e_3}$ component is in a PL form with the NEValue
\[
	\lambda _{3_{CL}}({\bf{e}}) =  - \left( {a{x^{*}_1}^2 + b{x^{*}_2}^2 + c({x_{s3}}^2 + {x_{s3}}{x_{m3}} + {x_{m3}}^2)} \right).
\] 
Since ${\lambda _{{3_{CL}}}}({\bf{e}}) \le 0$, the asymptotic exponential stability of $ e_{3}$ is guaranteed and then, the complete synchronization of all components will be obtained. The simulation results are shown in Figure~\ref{fig:synchronization}.

It is worth mentioning that for the parameter values in both the control and synchronization simulations, system~\eqref{eq:ghane} is chaotic.


\section{Concluding remarks}

The key idea of this paper was to use the PL form representation of nonlinear dynamical systems for the generation of chaotic behaviour. 
It is well known that the continual stretching and folding is the basic qualitative characteristic of a chaotic behaviour. This feature is essentially responsible for the local instability and global boundedness of chaotic trajectories. It has been shown that for a special class of nonlinear dynamical systems, the NEValues are indicators for the qualitative behaviour of the system. These qualitative indicators were applied to synthesize a particular form of continual stretching and folding behaviour in the state space of a 3-dimensional dynamical system. Numerical simulations verified the chaotic nature of the obtained system for a wide range of parameters. Chaotic dynamics arises through period doubling cascades of periodic attractors. Analysis by means of a Poincar\'e map suggests that the resulting chaotic attractors are of H\'enon-like type which means that they are the closure of an unstable manifold of a saddle periodic orbit. Due to symmetries the system also exhibits multi-stability which means that two different chaotic attractors coexist for the same parameter values.

In addition, we showed that by means of an eigen-structure based method the chaotic system can be easily both controlled and identically synchronized with another system through some nonlinear state feedback even if not all states are accessible. We tried to show that some efforts in nonlinear optimal control theory leading to the SDRE approach can be applied in another field of dynamical system theory. Currently, we are working on  the definition and control of nonlinear non-minimum phase system by the help of PL form representation and the results will be reported soon.


\end{document}